\documentclass[12pt]{amsart}

\usepackage[leqno]{amsmath}
\usepackage{amssymb}
\usepackage{amsxtra}
\usepackage{amscd}
\usepackage{hyperref}
\usepackage{graphicx, color}
\usepackage{dcpic}

\newtheorem{thm}{Theorem}

\newtheorem{defn}[thm]{Definition}

\newtheorem{lemma}[thm]{Lemma}
\newtheorem{cor}[thm]{Corollary}

\newcommand{\R}{\mathbb{R}}

\newcommand{\N}{\mathbb{N}}

\newcommand{\Z}{\mathbb{Z}}

\newcommand{\lra}{\longrightarrow}

\newcommand{\addFigure}[3][1.0]{\begin{figure}\includegraphics[scale=#1]{./#2}\caption{#3}\label{fig:#2}\end{figure}}

\newcommand{\addSection}[2]{\section{#2}\label{sec:#1}}

\newcommand{\addDefn}[2]{\begin{defn}{#2}\label{defn:#1}\end{defn}}
\newcommand{\refDefn}[1]{Section \ref{defn:#1}}

\title{On non-isotopic spanning surfaces for a class of arborescent knots}
\author{Lawrence Roberts}

\begin{document}
\begin{abstract}  
We use the methods of Hedden, Juhasz, and Sarkar to exhibit a set of arborescent knots that bound large numbers of non-isotopic minimal genus spanning surfaces. In particular, we describe a sequence of {\em prime} knots $K_{n}$ which will bound at least $2^{2n-1}$ non-isotopic minimal spanning surfaces of genus $n$.  
\end{abstract}
\maketitle

\addSection{surfaces}{Equivalence of surfaces}

\noindent In \cite{HJS}, M. Hedden, A. Juhasz, and S. Sarkar show that the knot $8_{9}$, with some orientation, is the oriented boundary of two distinct spanning surfaces, where we take two compact, connected, oriented surfaces $\Sigma_{1}, \Sigma_{2} \subset S^{3}$ to be equivalent if one is ambiently isotopic to the other. By taking the the $n$-fold connected sum of $8_{9}$ with itself, they obtain a sequence of knots $K_{n}, n \in \N$ which bound at least $n$ distinct spanning surfaces. The increasing number of surfaces arises, however, from the independent choices possible for each summand.\\
\ \\
\noindent In this paper, we will extend their technique of combining sutured Floer homology (\cite{JuhH}) with the Seifer form to improve this result to prime knots. Namely,

\begin{thm}\label{thm:main}
For each $n \in \N$, there exists an oriented, {\em prime} knot $K$ in $S^{3}$ which bounds at least $2^{2n-1}$ oriented minimal genus spanning surfaces, each with genus $n$. 
\end{thm}

\noindent Our examples will be arborescent knots which are also alternating. By a theorem of Menasco, a reduced, alternating, prime knot diagram represents a prime knot. $K$ being arborescent allows us to easily guarantee these conditions. The reader can find an example of the {\em prime} knots we construct in Figure \ref{fig:simpletree2}. That knot bounds at least $8$ surfaces of genus $2$.  \\
\ \\
\noindent Throughout this paper we will assume that:
\begin{enumerate}
\item $\Sigma$ is an oriented, compact, connected surface embedded in $S^{3}$, with a single boundary knot. 
\item $N(\Sigma)$ is a tubular neighborhood of $\Sigma$ equipped with a product structure by an orientation preserving diffeomorphism $\Sigma \times [-1,1] \hookrightarrow N(\Sigma) \subset S^{3}$
\item $\Sigma^{\pm}$ is the image of $\Sigma \times \{\pm 1\}$ under the inclusion of $N(\Sigma)$. Furthermore, if $a \subset \Sigma$ is an simple closed curve, then $a^{\pm}$ is the image of $a \times \{\pm 1\}$ in $\Sigma^{\pm}$. We will call these the positive/negative push-off(s) of $a$. 
\item $S^{3}_{\Sigma} = S^{3} \backslash \mathrm{Int} N(\Sigma)$
\end{enumerate}

\noindent We will use the following notions for the equivalence of oriented surfaces:

\addDefn{pairsmap}{A smooth map $f: (Y_{1}, \Sigma_{1}) \lra (Y_{2}, \Sigma_{2})$ is an orientation preserving diffeomorphism of pairs if $f:Y_{1} \lra Y_{2}$ is an orientation preserving diffeomorphism which induces an orientation preserving diffeomorphism from $\Sigma_{1} \subset Y_{1}$ to $\Sigma_{2} \subset Y_{2}$}

\addDefn{equivalence}{Two compact, connected, oriented, surfaces, $\Sigma_{1}$ and $\Sigma_{2}$, embedded in $S^{3}$ are {\em equivalent} if there is an orientation preserving isotopy $\Phi_{t}$ of $S^{3}$ with 1) $\Phi_{0}=\mathrm{Id}$ and 2) $\Phi_{1}$ restricting to a diffeomorphism from $\Sigma_{1}$ to $\Sigma_{2}$.}
\noindent After this paper first appeared, the authored learned that one can find examples among 2-bridge knots where there are $2^{2k-1}$ ``inequivalent'' incompressible Seifert surfaces for a knot of genus $k$. This follows from the work of Hatcher and Thurston in \cite{HaTh}. However, their notion of equivalence is different from that in this paper: in \cite{HaTh} two spanning surfaces are equivalent if they are isotopic in the complement of the knot. Similarly, Jessica Banks notes that M. Sakuma classified minimal genus Seifert surfaces for special arborescent links, which are very similar to the examples in this paper, again using isotopy in the complement of the link to provide the notion of equivalence. Using this classification she describes  examples of arborescent knots bounding $2^{2k-1}$ surfaces of minimal genus $k$ which are different up to isotopy {\em in the complement of the knot}, \cite{Bank}. In addition, this approach to classifying spanning surfaces use entirely different techniques that we will use in this paper.\\ 
\ \\
\noindent However, the notion of equivalence in definition \refDefn{equivalence} is stronger than that in \cite{HaTh}, \cite{Saku},  and \cite{Kaki}. Instead of requiring that the boundary of the surface $\Sigma$ lies on a fixed link $L$ -- in particular during an isotopy -- we instead use allow the boundary and the surface to be isotoped. Our notion of equivalence corresponds, therefore, more closely to that required for the isotopy classification of surfaces in $S^{3}$.\\
\ \\
\noindent As an illustration of the difference, J. R. Eisner provided examples of composite knots with infinitely many spanning surfaces, up to isotopies preserving the knot, and these examples give rise to Kakimizu complexes corresponding to the complex structure on $\R$ with $\Z$ as vertices, see section 3 in \cite{Kaki}. However, the different surfaces come from spinning one summand around the arc in the other that is used for the connect sum. This spinning can be undone if we allow the knot to move, so all these surfaces are equivalent under our definition.\\
\ \\
\noindent{\bf Acknowledgments:} The author would like to thank Allen Hatcher and Jessica Banks for some very helpful correspondence concerning the results of this paper and how they relate to previously known results. 

\addSection{trees}{Matched trees}

\subsection{Definitions}\ \\

\noindent Let $T$ be a finite, connected tree with vertices $V$ and edges $E$. In addition, suppose $T$ possesses this additional structure:
\begin{enumerate}
\item $T$ is bipartite:  $V$ is the disjoint union of $B$ and $W$ and any edge $e \in E$ has one vertex in $B$ and the other in $W$. 
\item $T$ has a matching: there is a bijection $m: B \rightarrow W$ such that for each $b \in B$, $b$ and $m(b)$ are the endpoints of an edge in $T$. (We will think of $m$ also as the subset of edges $\{b,m(b)\}$ in $E$)
\item $(T,m)$ is directed according to the convention that 1) $e \in m$ is oriented from its endpoint in $B$ to that in $W$, while 2) $e \not\in m$ will be directed oppositely, from $W$ to $B$.
\end{enumerate}
We will call this constellation of objects a matched tree $T$. The set of matched trees will be denoted $\mathcal{T}_{m}$. \\
\ \\
\noindent The simplest example of a matched tree is the tree with two vertices partitioned as $B = \{b\}$ and $W = \{w\}$ , equipped with the matching $m(b) = w$, and oriented accordingly. This particular tree will be denoted $T_{1}$.\\
\ \\
\addFigure[1]{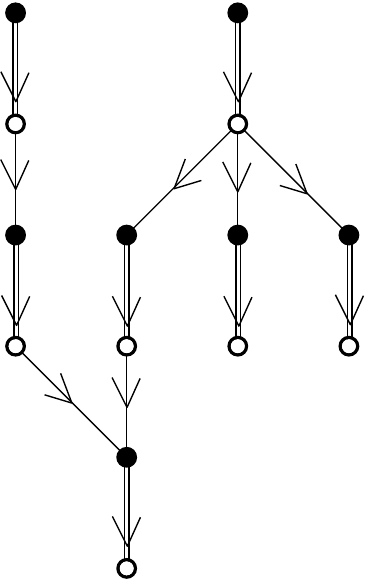}{A matched tree}

\noindent We will depict matched trees by drawing those edges in $m$ as double lines and edges in $E \backslash m$ as single lines. Furthermore, we identify those $v \in B$ by filled dots, while those $v \in W$ are identified by empty circles. For another example, drawn with this convention, see Figure \ref{fig:tree}.\\
\ \\
\noindent We can drop reference to the matching from $(T,m)$ since according to the following lemma there is a unique one, if any exist at all. 

\begin{lemma}
Any tree $T$ admits at most one matching.
\end{lemma} 

\noindent{\bf Proof:} If $m$ and $m'$ are distinct matchings on a bipartite graph $\Gamma$, then the symmetric difference $m \Delta m' \subset E$ consists of edges which form cycles in $\Gamma$. Consequently, when $\Gamma$ is a tree, there can be at most one matching.$\Diamond$\\
\ \\
\noindent Matched trees can also be described as those trees with the following properties

\begin{lemma}
Let $\mathcal{T}_{m}$ is the smallest set $\mathcal{T}$ of trees such that
\begin{enumerate}
\item $T_{1} \in \mathcal{T}$, and
\item $T \in \mathcal{T}$ if and only if there is a subtree $T'\in \mathcal{T}$ such that $V(T) = V(T') \cup \{l,v\}$, and $E(T) = E(T') \cup \{e_{1}, e_{2}\}$ with $e_{1}$ joining $v$ to a vertex in $T'$ and $e_{2}$ joining $l$ to $v$.
\end{enumerate}
\end{lemma}

\noindent We will see in the proof that the matched, bipartite structure on $T_{1}$ extends to all trees in this set. \\
\ \\
\noindent{\bf Proof:} First, we show that $\mathcal{T} \subset \mathcal{T}_{m}$. Since $T_{1}$ is in both, we suppose that any tree tree in $\mathcal{T}$ with $2n-2$ vertices is also in $\mathcal{T}_{m}$. We will then show that any tree in $\mathcal{T}$ with $2n$ vertices is also in $\mathcal{T}_{m}$, and draw the desired conclusion by induction. Suppose $T$ has $2n$ vertices, and that $T' \subset T$ is a subtree guaranteed by the second condition. Then $e_{1}$ joins $v$ to a vertex $v' \in T'$. First we extend the bipartite structure: if $v' \in B$ assign $v$ to $W$ and $l$ to $B$, otherwise assign $v$ to $B$ and $l$ to $W$. The matching on $T'$ uniquely extends by assigning $m(l) = v$ when $v' \in B$ (hence $l \in B$) and $m(v) = l$ when $v' \in W$. Finally, there is exists unique way to direct $T$ according to the orientation convention. Hence $T \in \mathcal{T}_{m}$ as well.  \\
\ \\
\noindent Now need we show that $\mathcal{T}_{m} \subset \mathcal{T}$. We again prove this by induction on the number of vertices in $T$. For $2$ vertices, the only connected, bipartite, matched tree is the starting tree $T_{1}$. Suppose $T$ has $2n$ vertices ($n \geq 2$) and that any matched tree in $\mathcal{T}_{m}$ with fewer than $2n$ vertices is in $\mathcal{T}$. Let $m$ be the matching for $T$. If we contract all the edges in $m$, we obtain a new tree $\Gamma$ with $\geq 2$ vertices.  $\Gamma$ has a leaf vertex $g$ which corresponds in $T$ to an edge $e_{2} \in m$ joining a leaf vertex $l$ to a bivalent vertex $v$. Let $e_{2}$ be the other edge at $v$. If we delete $e_{1},e_{2},v,l$ we obtain a tree $T'$ in $\mathcal{T}_{m}$ with $2n-2$ vertices. Thus $T' \in \mathcal{T}$, so $T \in \mathcal{T}$ as well.  $\Diamond$\\
\ \\
\noindent As a byproduct, this proof provides a way to iteratively construct all matched trees with a given number of vertices. 

\subsection{Partial Ordering}

\begin{defn}
Let $v, v'$ be vertices of a directed tree $T$. We will write $v \leq v'$ if there is a (possibly empty) directed path from $v'$ to $v$ in $T$. We will write  $v <_{1} v'$ if the directed path from $v'$ to $v$ contains exactly one edge.  
\end{defn}

\begin{lemma}
For any directed tree, $\leq$ is a partial order on the vertices of $T$. 
\end{lemma}

\noindent The orientation conventions for matched trees allow us to be more specific:

\begin{lemma}
When $T$ is a matched tree the terminal vertices for $(T,\leq)$ are all in $W$, while the initial vertices are all in $B$
\end{lemma}

\noindent For two vertices $v, v'$ in a directed tree the set $\big\{\,u \in V\,\big|\, v\leq u \leq v'\,\big\}$ is either empty (if $v \not\leq v'$) or consists of the vertices of the linearly ordered path from $v'$ to $v$. Furthermore, if $b \in B$ and $w \in W$ in a matched tree $(T,m)$, with $b \geq w$, then this path has a matching coming from the restriction of $m$. \\
\ \\
\noindent The bipartite structure on a matched tree allows us to refine some of these considerations:

\begin{defn}
For $w \in W$, let $B_{w} = \big\{\,b\in B\,\big|\, w \leq b\,\big\}$ be the vertices in $B$ that are above $w$. For $b \in B$,  let $W_{b} = \big\{\,w\in W\,\big|\, b \geq w\,\big\}$ be the vertices in $W$ that are below $b$. 
\end{defn}  

\noindent For a matched tree, the partial order on $B$ and $W$ is reflected in these sets alone.

\begin{lemma}
Let $T$ be a bipartite tree oriented by $m$. For $b, b' \in B$, $b \leq b'$ if and only if $W_{b} \subset W_{b'}$. Likewise, for $w, w' \in W$, $w \leq w'$ if and only if $B_{w} \supset B_{w'}$.  
\end{lemma}


\section{The surfaces}

\noindent Let $T$ be a matched tree as in the previous section. Equip $T$ with 
\begin{enumerate}
\item a {\em framing} map $f : V \rightarrow \Z$, and
\item a {\em plumbing} map $\epsilon : E \rightarrow \{-1,+1\}$
\end{enumerate} 
Given $(T,f,\epsilon)$, we can define an oriented surface in $S^{3}$. To do this recall

\addDefn{murasugi}{An oriented surface $\Sigma \subset S^{3}$ is the Murasugi sum of the oriented surfaces $\Sigma_{1}$ and a surface $\Sigma_{2}$ if there is a decomposition $S^{3} = B_{1} \cup_{S} B_{2}$ where $B_{1}$ and $B_{2}$ are closed balls, such that
\begin{enumerate}
\item For $i=1,2$, $B_{i} \cap \Sigma \subset S^{3}$ is isotopic to $\Sigma_{i}$, while preserving the orientations on the surfaces,
\item $S \cap \Sigma$ is a disc $D$ such that $\partial D \cap \partial \Sigma$ consists of $2n > 2$ points on is boundary.
\end{enumerate}
}
\addDefn{pos}{The Murasugi sum is positive for $\Sigma_{i}$ if the orientation for $S$ inherited from $S \cap \Sigma$ agrees with its orientation as the boundary of $B_{i}$. It is negative, if it disagrees.} 

\addFigure[0.75]{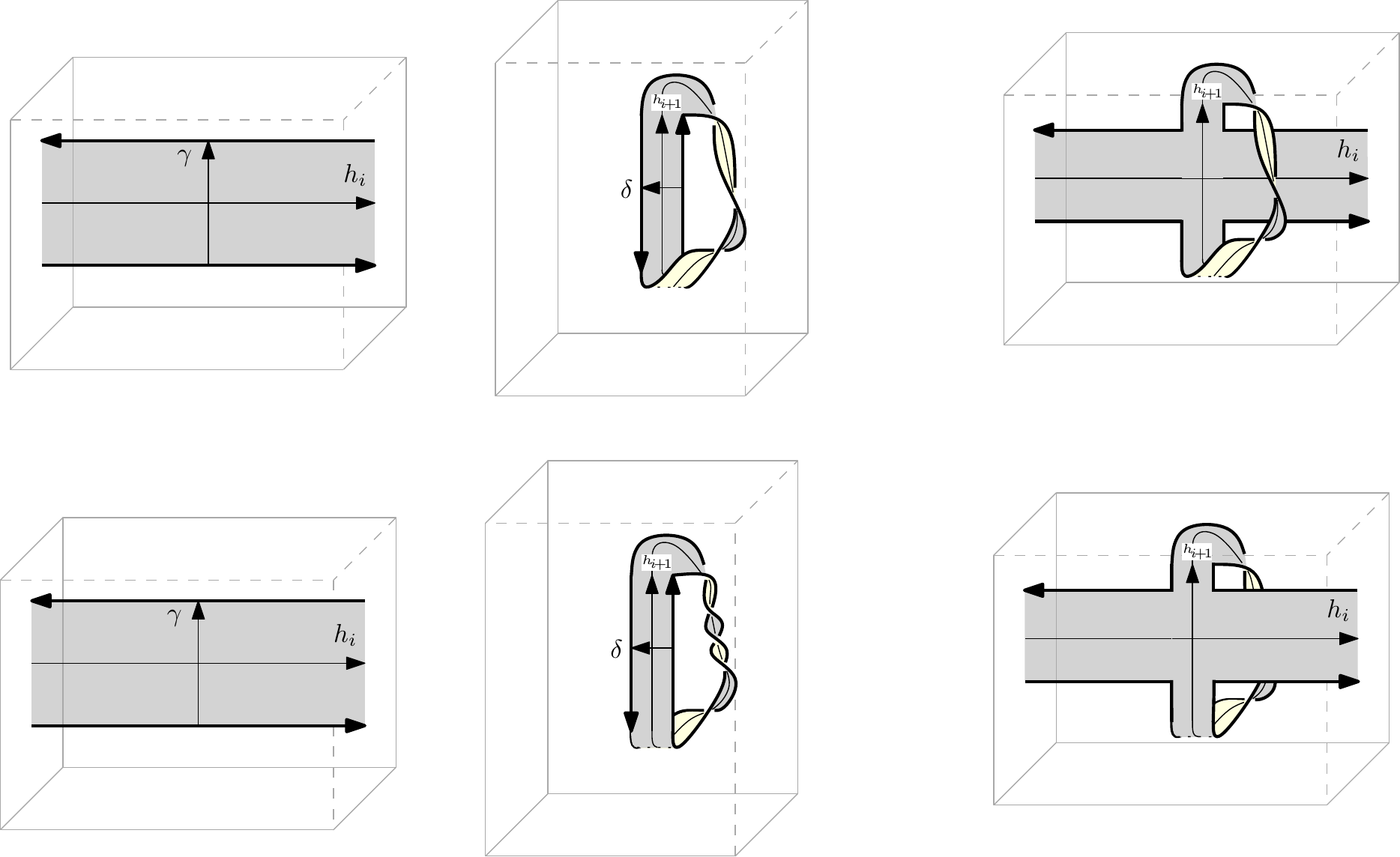}{The top digram is a positive plumbing of the annulus, using the notation below. The bottom is a negative plumbing. Note that the annuli have the same number of full twists in each. The boundaries of the two surfaces on the right will represent the same link class.}

\noindent We provide two local models for the $n=2$ case of the Murasugi sum. Two diagrams illustrate these models in Figure \ref{fig:plumbings} Suppose that $\Sigma \subset S^{3}$ is a compact, oriented surface with boundary. We assume we are given an oriented circle $h$ in $\Sigma$, and an oriented arc $\gamma$, properly embedded in $\Sigma$, with $h \pitchfork \gamma = +1$, algebraically and geometrically. Let $N$ be a neighborhood of $\gamma$ in $S^{3}$ equipped with coordinates such that
\begin{enumerate}
\item  N is homeomorphic to $[-1,1]^{3}$,
\item $\gamma$ is  $\{0\} \times [-1/2,1/2] \times \{0\}$,
\item $h \cap N$ is $[-1,1] \times \{(0,0)\}$
\item $\Sigma \cap N$ is $[-1,1] \times [-1/2,1/2] \times \{0\}$ with $\partial \Sigma \cap N = \big([-1,1] \times \{(-1/2,0)\}\big) \cup -\big([-1,1] \times\{(1/2,0)\}\big)$. 
\end{enumerate} 
and all identifications preserve the relevant orientations. \\
\ \\
\noindent Let $\Sigma' \subset S^{3}$ be another surface equipped with a closed curve $C$ and an $\delta$ from the boundary to the boundary, so that $C \pitchfork \delta = +1$, algebraically and geometrically. Let $\epsilon \in \{1,-1\}$. Suppose $M \subset S^{3}$ is homeomorphic to $[-1,1]^{3}$ such that
\begin{enumerate}
\item $\delta$ equals $-\big([-1/2,1/2] \times \{(0,0)\}\big)$. 
\item $C \cap [-1,1]^2 \times\{0\} = \{0\} \times [-3/4,3/4] \times \{0\}$. 
\item $\Sigma' \cap [-1,1]^2 \times\{0\}$ is $[-1/2,1/2] \times [-3/4,3/4] \times \{0\}$ with $[-1,1]^2 \times\{0\} \cap \partial A = -\big(\{-1/2\} \times [-1,1] \times \{0\}\big) \cup \big(\{1/2\} \times [-1,1] \times \{0\}\big)$
\item If $\epsilon = +1$ then $\Sigma' \subset$ $[-1,1]^{2} \times [0,1]$, whereas if $\epsilon = -1$ then $\Sigma' \subset$ $[-1,1]^{2} \times [-1,0]$
\end{enumerate}

\noindent Then $\Sigma \ast_{\gamma, \delta, \epsilon} \Sigma'$ is the union of $\Sigma$ with the image of $\Sigma'$ under the coordinate identification of $N$ and $M$. This surface is a Murasugi sum of $\Sigma$ and $\Sigma'$ with decomposing sphere equal to $S_{+1} = \partial\big([-1,1]^{2} \times [0,1]\big)$ when $\epsilon = +1$, and $S_{-1} = \partial\big([-1,1]^{2} \times [-1,0]\big)$ when $\epsilon = -1$. The sum is thus positive for $\Sigma$ when $\epsilon = +1$, and negative for $\Sigma$ when $\epsilon = -1$. \\
\ \\
\noindent Now let $(T,f,\epsilon)$ be a matched tree equipped with vertex and edge labellings (which we will shorten to $T$ from here on). For each vertex $v$ of $T$ let $A_{v}$ be an annulus in $S^{3}$, with {\em unknotted}, oriented, core circle $h_{v}$, making $f(v)$ right handed twists around $h_{v}$ relative to the framing provided by a disk $h_{v}$ bounds. If $v$ is the endpoint for the edges $e_{1}, \ldots, e_{n}$ pick simple arcs $\gamma_{e_{1}}, \dots, \gamma_{e_{n}}$ running between the two boundary components, with $h_{v} \pitchfork \gamma_{i} = +1$. As the constructions will not depend on these arcs, we will suppress them from the notation. \\
\ \\
\noindent To construct $\Sigma_{(T,f,\epsilon)}$ we Murasugi sum  all the annuli $A_{v}$ by using the arcs $\gamma_{e} \subset A_{v}$ and $\gamma_{e}' \subset A_{v'}$ corresponding to edge $e$ between $v$ and $v'$. To fully specify $\Sigma_{T}$ we need to know whether to use the positive or negative sum for each edge. For edge $e$ from $v$ to $v'$, with $v < v'$, we use the local model from $\epsilon(e)$ with $\Sigma$ being the surface containing $A_{v}$ and $\Sigma'$ being the surface containing $A_{v'}$ for the larger vertex $v'$.  

\begin{lemma}
For a surface $\Sigma_{T}$ as constructed above, the following are true
\begin{enumerate}
\item $S^{3}_{\Sigma_{T}}$ is homeomorphic to a handlebody.
\item $\partial \Sigma_{T}=K_{T}$ is a knot
\item The knot type of $K_{T}$ does not depend on the map $\epsilon$, when $T$, $f$, and the arcs $\gamma_{e}$ are all fixed.
\end{enumerate}  
\end{lemma}

\noindent{\bf Proof:} All three are proved by induction. First, they are true for triples $(T_{1},f_{1},\epsilon_{1})$ for any $f_{1}$ and $\epsilon_{1}$ can be verified by directly checking. Now suppose that for any $T' \in \mathcal{T}_{m}$ with fewer than $2n$ vertices, all three properties are true. Let $T \in \mathcal{T}_{m}$ have $2n$ vertices. Then there are two edges $e_{1}$ and $e_{2}$ with $e_{2} \in m$, joining $l$ to $v$. And $e_{1}$ joining $v$ to $v'$, such that $T$ minus $e_{1},e_{2},v,l$ is in $\mathcal{T}_{m}$. By construction $e_{1}$ corresponds to a sphere $S \subset S^{3}$ which effects the sum of $\Sigma_{T'}$ to the surface specified by $v$, $l$, and $e_{2}$, which independently form a copy of $T_{1}$ in one of the balls $B$ with boundary $S$. The complement of this surface in $B$ consists of two one handles attached to a thickening of $S$, and thus when glued to the complement of $\Sigma_{T'}$ in the other ball, will still be a handlebody (since $B$ can be thought of as a portion of the neighborhood of $\Sigma_{T'}$. Since $\partial \Sigma_{T_{1}}$ and $\partial \Sigma_{T'}$ are both knots, the sum will also be a knot. Finally, the link that is the boundary of the sum of the annulus for $v$ to the surface $\Sigma_{T'}$ will not depend on the choice of $\epsilon$ from the local models above, and likewise for the sum of $l$ to $\Sigma_{T' \cup \{v\}}$. The proposition is thus proved. $\Diamond$\\
\ \\
\noindent Note that $K_{T}$ is an arborescent knot, by construction. There is an example of a tree $T$ and the corresponding knot $K_{T}$ in Figure \ref{fig:simpletree2}.

\addFigure[1.0]{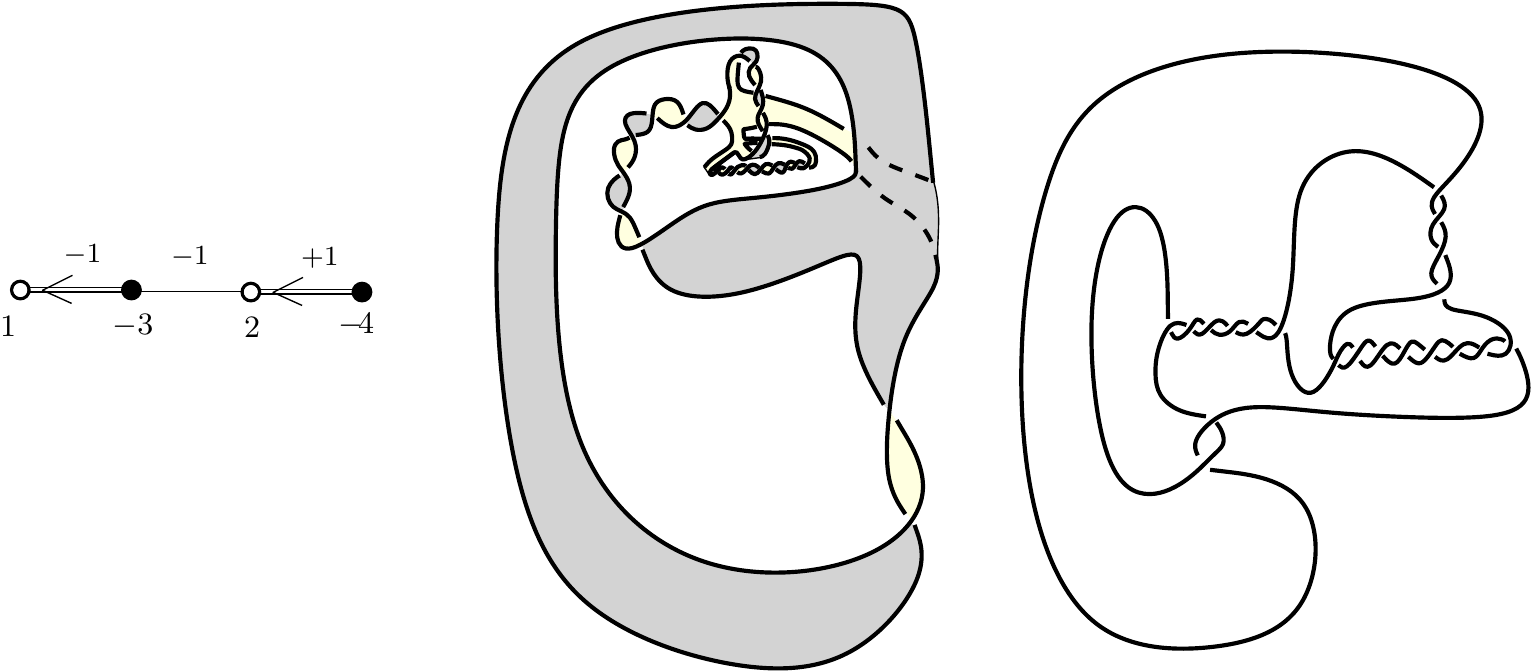}{}


\section{Seifert Forms}

\noindent For each $\Sigma_{T,\epsilon}$ there is a Seifert form $\theta_{T,\epsilon} : H_{1}(\Sigma;\Z) \otimes H_{1}(\Sigma; \Z) \rightarrow \Z$ defined by  $$\theta_{T,\epsilon}([a],[b]) = lk(a,b^{+})$$ 
\noindent This form is invariant under orientation preserving diffeomorphisms of pairs $f: (S^{3},\Sigma_{1}) \lra (S^{3}, \Sigma_{2})$.\\
\ \\
\noindent Indeed, such a diffeomorphism of pairs maps the product structure of $N(\Sigma_{1})$ to one for $N(\Sigma_{2})$. Consequently, $f(a^{+}) = (f(a))^{+}$ for any cycle $a \subset \Sigma_{1}$. As it will likewise map surfaces with boundary $a$ or $a^{+}$ we see that 
$$\theta_{\Sigma_{2}}(f_{\ast}(\alpha_{1}), f_{\ast}(\alpha_{2})) = \theta_{\Sigma_{1}}(\alpha_{1},\alpha_{2})$$ where $f_{\ast}$ is the map induced on the homology of $\Sigma_{1}$. \\
\ \\
\noindent To find the Seifert form for $\Sigma_{T,\epsilon}$ we will use the local models for positive/negative Murasugi sum. Note that the construction of $\Sigma_{T,\epsilon}$ also provides a basis $\{h_{v}|v \in V\}$ for $H_{1}(\Sigma_{T,\epsilon})$ with the property that 
$$
h_{v} \pitchfork h_{v'} = \left\{\begin{array}{cl} +1 & v <_{1} v' \\ 0 & \mathrm{otherwise} \end{array}\right.
$$

\begin{lemma}
The Seifert form for $\Sigma_{T,\epsilon}$ is the form $\theta_{T,\epsilon}$ where
\begin{enumerate}
\item for each $v \in V$ let $\theta_{T}(h_{v}, h_{v}) = f(v)$,
\item for $v, v' \in V$ and $\{v,v'\} \not\in E$ let $\theta_{T}(h_{v}, h_{v'}) = \theta_{T}(h_{v'}, h_{v}) = 0$
\item for $v, v' \in V$ with $v <_{1} v'$, if
\begin{enumerate}
\item $\epsilon(\{v,v'\}) = +1$, then $\theta_{T}(h_{v}, h_{v'}) = 0$ while $\theta_{T}(h_{v'}, h_{v}) = 1$
\item $\epsilon(\{v,v'\}) = -1$, then $\theta_{T}(h_{v}, h_{v'}) = -1$ while $\theta_{T}(h_{v'}, h_{v}) = 0$
\end{enumerate}
\end{enumerate}
\end{lemma}

\noindent{\bf Proof:} If $v$ and $v'$ are not joined by an edge in $T$ then $lk(h_{v},h_{v'}^{+}) = lk(h_{v'},h_{v}^{+}) = 0$, since the annuli $A_{v}$ and $A_{v'}$ are located in disjoint balls. Furthermore, for $lk(h_{v}, h_{v}^{+})$ we need only consider the situation for the annulus $A_{v}$, where it is clear that we recover the framing $f(v)$. For $v <_{1} v'$ we will use the local models. If $\epsilon = +1$ then the surface for $v'$ lies in a ball on the positive side of the image of $A_{v}$. Consequently, $h_{v'}^{+}$ lies entirely within this ball, and bounds a disc there. So $lk(h_{v},h_{v'}^{+})= 0$. When $\epsilon = -1$, the surface lies on the negative side. Thus $h_{v'}^{+}$ loops around $h_{v}$ in the same manner as $(cos(\theta), 0, sin(\theta))$ loops around the $y$-axis in $\R^{3}$, namely in a left-handed manner. Thus, when $\epsilon = -1$  we have $lk(h_{v},h_{v'}^{+})= -1$. By the properties of the Seifert form $\theta_{\Sigma}(h_{v}, h_{v'}) - \theta(h_{v'}, h_{v}) = - h_{v} \pitchfork h_{v'} = +1$. Thus when $\epsilon = +1$, $\theta(h_{v'}, h_{v}) = 1$, but when $\epsilon = -1$,  $\theta(h_{v'}, h_{v}) = -1 + 1 = 0$. This is exactly the form $\theta_{T}$ defined previously. $\Diamond$.  \\
\ \\
\noindent Since $K_{T} = \partial \Sigma_{T}$ is a knot the groups $H_{1}(\Sigma_{T};\Z)$ and $H_{1}(S^{3}_{\Sigma_{T}};\Z)$ are isomorphic (as in \cite{HJS} or \cite{Kauf}). In fact, there is an explicit isomorphism provided in \cite{Kauf}: $\Phi_{\Sigma} : H_{1}(\Sigma;\Z) \lra H_{1}(S^{3}_{\Sigma};\Z)$ where 
$$\Phi_{\Sigma}([a]) = [a^{+}] - [a^{-}]$$
This isomorphism has the following naturality:
if  $f\!: (S^{3},\Sigma_{1}) \lra (S^{3}, \Sigma_{2})$ is an orientation preserving diffeomorphism of pairs, where the surfaces have a single boundary component, then the following commutes
$$
\begin{CD}
H_{1}(\Sigma_{1};\Z) @>\Phi_{\Sigma_{1}}>> H_{1}(Y_{1,\Sigma_{1}};\Z)\\
@V(f|_{\Sigma_{1}})_{\ast}VV               @VV(f|_{Y_{1,\Sigma_{1}}})_{\ast}V\\
H_{1}(\Sigma_{2};\Z) @>\Phi_{\Sigma_{2}}>> H_{1}(Y_{2,\Sigma_{2}};\Z)\\
\end{CD}
$$
since $f$ preserves the positive and negative push-offs. Consequently, $$\overline{\theta}_{\Sigma}([\alpha],[\beta]) = \theta_{\Sigma}(\Phi_{\Sigma}^{-1}([\alpha]), \Phi_{\Sigma}^{-1}([\beta]))$$  is a pairing on $H_{1}(S^{3}_{\Sigma};\Z)$ which will  be preserved by the restriction to $S^{3}_{\Sigma_{1}}$ of an orientation preserving diffeomorphism of pairs, in the same manner as $\theta_{\Sigma}$. \\
\ \\
\noindent We now compute $\overline{\theta}_{T,\epsilon}$. First we describe the map $\Phi$.

\begin{lemma}
Given a matched tree $T$ with framing $f$ and plumbing $\epsilon$, let $\{c_{v}|v \in V\}$ be the basis for $H_{1}(S^{3}_{\Sigma_{T,\epsilon}})$ given by meridians to the annuli $A_{v}$, and oriented so that $\mathrm{lk}(c_{v}, h_{v}) = +1$. Using these bases, the isomorphism $\Phi_{\Sigma_{T,\epsilon}}$ is defined by
$$
h_{v} \longrightarrow \sum_{v <_{1} v'} c_{v'} - \sum_{v' <_{1} v} c_{v'}
$$
\end{lemma}

\noindent{\bf Proof:} Let $h_{v}$ be a basis element for $H_{1}(\Sigma_{T})$.  We saw above that $h_{v} \pitchfork h_{v'} = +1$ if $v <_{1} v'$ and $h_{v} \pitchfork h_{v''} = -1$ if $v'' <_{1} v$. Then $h_{v}^{+} - h_{v}^{-}$ bounds an annulus $C$ which intersects each $h_{v'}$ and $h_{v''}$. This annulus has orientation given by a positive normal to $A_{v}$, followed by a tangent to $h_{v}$. Consequently, for $h_{v} \pitchfork h_{v'} = +1$, the triple $N, h_{v}, h_{v'}$ is a positively oriented basis. Therefore, $h_{v'}$  intersects $C$ positively. Similarly, for $v'' <_{1} v$ the intersection is negative. Thus $C$ minus open neighborhoods of the points of intersection provides a homology relation showing that $\Phi_{\Sigma_{T}}([h_{v}])$ is homologous to the the sum of the meridians around $h_{v'}$ for $v <_{1} v'$ and minus the meridians around $h_{v''}$ for $v'' <_{1} v$. These meridians are precisely $c_{v}$, so 
$$
\Phi_{\Sigma_{T}}(h_{v}) = \sum_{v <_{1} v'} c_{v'} - \sum_{v' <_{1} v} c_{v'}
$$ 
as required. $\Diamond$\\
\ \\
\noindent To compute $\overline{\theta}_{T,\epsilon}$,  however, we need $\Phi^{-1}$.

\begin{lemma}
In the bases described above, $\Phi^{-1}$ is generated by
\begin{align*}
c_{b} &\longrightarrow \sum_{w \in W_{b}} h_{w} \mathrm{\ \ for\ }b\in B\\
c_{w} &\longrightarrow -\sum_{b \in B_{w}} h_{b} \mathrm{\ \ for\ }w\in W\\
\end{align*}
\end{lemma}
 
\noindent{\bf Proof:} Let $b \in B$ and $w = m(b)$. Then $h_{m(b)} \rightarrow c_{b} - \sum_{v <_{1} m(b)} c_{v}$. Since $T$ is bipartite, each $v$ with $v <_{1} m(b)$ will also be in $B$. In particular, $v < b$ so $W_{v} \subset W_{b}$. If we rewrite this as $c_{b} \rightarrow h_{m(b)} + \sum_{v <_{1} m(b)} c_{v}$ and repeat with each of the $v$'s we will obtain $c_{b} \rightarrow \sum_{w \in W_{b}} h_{w}$. For $w \in W$ we know that $w=m(b)$ for some $b$ and that $w$ is the only vertex satisfying $b >_{1} v$. Consequently, $h_{b} \rightarrow \sum_{b <_{1} v} c_{v} - c_{w}$. Rearranging we obtain $c_{w} \rightarrow \sum_{b <_{1} v} c_{v} - h_{b}$. Each $v$ with $b <_{1} v$ is in $W$, so we can repeat with each term in the first sum to obtain $c_{w} \rightarrow -\sum_{b \in B_{w}} h_{b}$. $\Diamond$. \\
\ \\
\noindent From this description we can calculate $\theta_{T}(\Phi^{-1}(c_{u}), \Phi^{-1}(c_{v}))$: for $b, b' \in B$ 
$$
\overline{\theta}_{T}(c_{b},c_{b'}) = \sum_{w \in W_{b} \cap W_{b'}} f(w)
$$
while for $w,w' \in W$
$$
\overline{\theta}_{T}(c_{w},c_{w'}) = \sum_{b \in B_{w} \cap B_{w'}} f(b)
$$
Furthermore, for $v, v'$ let $\gamma$ be the directed path from $v'$ to $v$ (or be empty if there isn't one). Let $N_{\pm, m}(\gamma)$ be the number of $\epsilon = \pm 1$ edges on $\gamma$ which are also in $m$, and $N_{\pm,\sim m}$ be the number that are not in $m$. Then
for $b \in B$ and $w \in W$, let $\gamma$ be the directed path from $b$ to $w$, or else empty if no such path exists. We have\\
$$
\overline{\theta}_{T}(c_{w},c_{b}) = N_{-1,\sim m}(\gamma) - N_{+1, m}(\gamma)
$$
$$
\overline{\theta}_{T}(c_{b},c_{w}) = N_{-1, m}(\gamma) - N_{+1,\sim m}(\gamma)
$$
where it is to be understood that if $w > b$ then both of these are $0$. \\
\ \\
\noindent We illustrate the latter calculation. Suppose $w < b$. Starting from $w$ we will enumerate the vertices up the path until we arrive at $b$. Thus we get $w=v_{1}, \ldots, v_{2g} = b$ for some $g$. Since $\theta_{T}$ is zero on basis elements, unless they correspond to consecutive vertices, we can see that only the edges in this path contribute to $\overline{\theta}_{T}(c_{w},c_{b})$, thus\\
\begin{align*}
\overline{\theta}_{T}(c_{w},c_{b}) &= - \big(\theta_{T}(h_{2},h_{1}) + \theta_{T}(h_{2}, h_{3}) + \theta_{T}(h_{4},h_{3}) + \theta_{T}(h_{4}, h_{5}) +  \ldots \\
&\hspace{0.5in}  + \theta_{T}(h_{2g-2}, h_{2g-1}) + \theta_{T}(h_{2g},h_{2g-1})\big)\\  
\end{align*}
The edge from $v_{2i}$ to $v_{2i-1}$ is in $m$, while that from $v_{2i+1}$ to $v_{2i}$ is not. Consequently, those edge in $m$ contribute when they evaluate to $+1$ under $\epsilon$ since the even number comes first. In this case, the value of $\theta_{T}$ will be $+1$ for this edge, so with the minus sign in front the edges in $m$ contribute $-N_{+1,m}(\gamma)$. The edges not in $m$ contribute a $-1$ when $\epsilon = -1$, so altogether they contribute $-(-N_{-1,\sim m}(\gamma)) = N_{-1,\sim m}(\gamma)$.

\section{Sutured manifolds and Floer homology}

\noindent $S^{3}_{\Sigma} = S^{3} \backslash \mathrm{Int}N(\Sigma)$ is a sutured manifold in the sense of Gabai (see \cite{HJS}). Indeed from the product structure on $N(\Sigma)$ we have $\partial S^{3}_{\Sigma} =$  $-\Sigma^{+} \cup$, $\Sigma^{-} \cup$ $\partial \Sigma \times [-1,1]$. Thus $(S^{3}_{\Sigma}, \Sigma^{+}, \Sigma^{-})$ is a balanced sutured manifold ( \cite{HJS}) with suture equal to the annulus $\partial \Sigma \times [-1,1]$ and $R(\gamma) = \Sigma^{+} \cup \Sigma^{-}$\footnote{However, in the notation of \cite{HJS} $R_{+}(\gamma) = \Sigma^{-}$ and $R_{-}(\gamma) = \Sigma^{+}$ since the oriented normal to $\Sigma^{+}$ points into $S^{3}_{\Sigma}$}.  \\
\ \\
\noindent An orientation preserving diffeomorphism $f: M_{1} \lra M_{2}$ between sutured manifolds $(M_{1},\gamma_{1})$ and  $(M_{2}, \gamma_{2})$ is called an identification of sutured manifolds if it restricts to $R(\gamma_{1}) = (\partial M_{1}) \backslash \gamma_{1}$ in a manner preserving the orientation assigned to each component by the sutures. Applied to our setting, $f: S^{3}_{\Sigma_{1}} \rightarrow S^{3}_{\Sigma_{2}}$ must restrict to each of $\Sigma_{1}^{\pm}$ as an orientation preserving diffeomorphism with $\Sigma_{2}^{\pm}$. Consequently, orientation preserving diffeomorphisms of pairs $(S^3,\Sigma)$ induce identifications of the corresponding sutured manifolds $(S^3_{\Sigma}, \Sigma^{+}, \Sigma^{-})$. \\
\ \\
\noindent A. Juhasz, \cite{JuhH}, associates a $\Z/2\Z$ vector space $\mathrm{SFH}(M, \gamma; \mathfrak{s})$ to any balanced sutured manifold $(M,\gamma)$ and a choice of (relative) $\mathrm{Spin^{c}}$ structure $\mathfrak{s} \in \mathrm{Spin}^{c}(M,\gamma)$. We will not review the definition of sutured Floer homology here. Instead we list several properties of this invariant.\\
\ \\
\noindent First, recall that there is an action of $H_{1}(M;\Z)$ on $\mathrm{Spin}^{c}(M,\gamma)$ such that for any $\mathfrak{s}_{1}, \mathfrak{s}_{2} \in \mathrm{Spin}^{c}(M,\gamma)$ there exists a unique $h \in H_{1}(M;\Z)$ with $h \cdot \mathfrak{s}_{1} = \mathfrak{s}_{2}$. We will denote this by $h=\mathfrak{s}_{2} - \mathfrak{s}_{1}$. Then  

\begin{enumerate} 
\item $\mathrm{SFH}(M, \gamma; \mathfrak{s}) \not\cong \{0\}$ for finitely many $\mathfrak{s} \in \mathrm{Spin}(M,\gamma)$.
\item \label{list2} Any diffeomorphism $f: (M_{1},\gamma_{1}) \rightarrow (M_{2}, \gamma_{2})$ of sutured manifolds induces a bijection $F^{\mathrm{Sp}}: \mathrm{Spin}^{c}(M_{1},\gamma_{1}) \rightarrow \mathrm{Spin}^{c}(M_{2},\gamma_{2})$ such that $F^{\mathrm{Sp}}(\alpha \cdot \mathfrak{s}) = f_{\ast}(\alpha) \cdot F^{\mathrm{Sp}}(\mathfrak{s})$, where $\alpha \in H_{1}(M_{1};\Z)$, 
\item For a diffeomorphism $f$ as in item (\ref{list2}) and any $\mathfrak{s} \in \mathrm{Spin}^{c}(M_{1},\gamma_{1})$ there is an isomorphism $\mathrm{SFH}(M_{1},\gamma_{1}; \mathfrak{s}) \stackrel{f_{\mathfrak{s}}}{\lra} \mathrm{SFH}(M_{2},\gamma_{2}; F^{\mathrm{Sp}}(\mathfrak{s}))$
\end{enumerate}
\noindent We note that $\mathrm{SFH}(M, \gamma; \mathfrak{s})$ also comes with an invariant relative grading, but we will not need this in the sequel.\\
\ \\
\noindent We will now describe $\mathrm{SFH}(S^{3}_{\Sigma_{T}})$ for a surface from a labeled, matched tree. Since these surfaces arise from repeated Murasugi sums we will use proposition  8.6 of \cite{JuhS}, which shows that
$$
\mathrm{SFH}(S^{3}_{\Sigma}) \approx \mathrm{SFH}(S^{3}_{\Sigma_{1}}) \otimes \mathrm{SFH}(S^{3}_{\Sigma_{2}})
$$
when $\Sigma$ is the Murasugi sum of $\Sigma_{1}$ and $\Sigma_{2}$. \\
\ \\
\noindent In addition, when $\Sigma$ is the Murasugi sum of $\Sigma_{1}$ and $\Sigma_{2}$,  $\mathrm{Spin}^{c}(S^{3}_{\Sigma}) = \mathrm{Spin}^{c}(S^{3}_{\Sigma_{1}}) \times \mathrm{Spin}^{c}(S^{3}_{\Sigma_{2}})$, and this is a torsor of $H_{1}(S^{3}_{\Sigma};\Z) \cong H_{1}(S^{3}_{\Sigma_{1}}; \Z) \oplus H_{1}(S^{3}_{\Sigma_{2}};\Z)$, with  $H_{1}(S^{3}_{\Sigma_{2}};\Z)$ acting solely on the $\mathrm{Spin}^{c}(S^{3}_{\Sigma_{2}})$ factor, and $H_{1}(S^{3}_{\Sigma_{1}}; \Z)$ acting solely on the $\mathrm{Spin}^{c}(S^{3}_{\Sigma_{1}})$ factor.  \\
\ \\
\noindent To apply this to the sums defining $\Sigma_{T}$ we also need the sutured Floer homology for the complement of an unknotted, twisted annulus in $S^{3}$. Fortunately, A. Juhasz computed this in \cite{JuhP}. \\
\ \\
\noindent  Let $A_{p} \subset S^{3}$ be the unknotted band with $p \in \Z$ right handed full twists. Then $N(A_{p}) \cong S^{1} \times D^{2}$ and $S^{3} \backslash \mathrm{Int} N(A_{p}) \cong D^{2} \times S^{1}$. The two components of $\partial A_{p}$ give two curves on $\partial N(A_{p})$ which go around the longitude once, and $p$ times around the meridian, of the unknotted circle forming the core of $A_{p}$. Thinking of the the boundary torus as the boundary of $S^{3} \backslash \mathrm{Int} N(A_{p})$ switches the meridian and longitude. Thus, in $S^{3} \backslash \mathrm{Int} N(A_{p}) \cong D^{2} \times S^{1}$, the two curves provide sutures which wind around the boundary of $D^{2}$ once and $p$ times around the $S^{1}$ factor.\\
\ \\
\noindent Thus the sutured manifold $S^{3} \backslash \mathrm{Int} N(A_{p})$ will be diffeomorphic to $T(p,1;2)$, where, following A. Juhasz, \cite{JuhP}, $T(p,q;n)$ is the sutured manifold $(S^{1} \times D^{2},\gamma)$ with suture $\gamma$ consisting of $n=2k+2$ parallel disjoint curves in the homology class $(p,q) \in H_{1}(S^{1} \times S^{1};\Z)$. Proposition 9.1 of \cite{JuhP} computes the sutured Floer homology of $T(p,q;n)$ in general. For $T(p,1;2)$ this proposition states that there is an identification $\mathrm{Spin}^{c}(T(p,1;2)) \leftrightarrow \Z$ such that
$$
\mathrm{SFH}(T(p,1;2),i) \cong \left\{\begin{array}{cl} \Z & 0 \leq i < p \\
0 & \mathrm{otherwise}\\
\end{array} \right.
$$
Furthermore, the difference of two structures in $\mathrm{Spin}^{c}(T(p,1;2))$ is a multiple of the meridian $c$ of the twisted band. \\
\ \\
\noindent Repeatedly applying these propositions to the Murasugi sums defining $\Sigma_{T}$,  we see that 
$$
Q_{T}= \big\{\mathfrak{s} \in \mathrm{Spin}^{c}(S^{3}_{\Sigma_{B}}) \big|  \mathrm{SFH}(S^{3}_{\Sigma_{B}}) \not\cong 0\big\} \cong \prod_{v \in V} [0,|f(v)|-1]
$$
The product structure reflects the action of the meridional basis $\{c_{v}|v\in V\}$. Namely, there is a $\mathrm{Spin}^{c}$ structure $\mathfrak{s}$ so that all the other structures with non-zero sutured Floer homology can be found by taking $\mathfrak{s} + \sum a_{v} c_{v}$ with $a_{v}$ in the interval for $[0,|f(v)|-1]$. In particular, the lengths of the intervals come from the identifications in the calculation of $\mathrm{SFH}(T(p,1;2))$. \\

\section{Symmetries}

\noindent Finally we prove

\begin{thm}
Let $T$ be a matched tree with $2n$ vertices, equipped with an framing map $f: V \rightarrow \Z\backslash\{0, \pm 1\}$ such that $v \rightarrow |f(v)|$ is injective. Suppose further that gluing arcs are fixed for each annulus $A_{v}$, sufficient to implement the construction of $\Sigma_{T, \epsilon}$ for any edge labeling $\epsilon: E \rightarrow \{\pm 1\}$. Then, using the same choice of arcs, $\Sigma_{T,\epsilon}$ is not equivalent to $\Sigma_{T, \epsilon'}$ unless $\epsilon = \epsilon'$.  
\end{thm}

\noindent {\bf Note:} We disallow $f(v)$ to have absolute value $1$ so that $Q_{\Sigma_{T}}$ will have dimension $2n$. 

\begin{cor}
For $T,f$ as above, the knot $K_{T}$ is the boundary of $2^{2n-1}$ distinct embedded surfaces of genus $n$. 
\end{cor}

\noindent {\bf Proof:} We show that there is no orientation preserving diffeomorphism of pairs $f:(S^{3},\Sigma_{T,\epsilon}) \lra (S^{3},\Sigma_{T,\epsilon'})$. To do this, we start by describing the possible maps $f_{\ast}: H_{1}(S^{3}_{\Sigma_{T,\epsilon}};\Z) \lra H_{1}(S^{3}_{\Sigma_{T,\epsilon'}};\Z)$ in terms of the meridional bases $\{c_{v}\}$ and $\{c'_{v}\}$. These maps are severally constrained by the properties of sutured Floer homology and the injectivity of $f$. \\
\ \\
\noindent The structure of $Q_{T,\epsilon}$ does not depend on $\epsilon$. Consequently, a map  $f: (S^{3},\Sigma_{T,\epsilon}) \lra (S^{3},\Sigma_{T,\epsilon'})$ must induce a map $F^{Sp}$ which takes the supporting parallelepiped for $S^{3}_{\Sigma_{T,\epsilon}}$ to that for $S^{3}_{\Sigma_{T,\epsilon'}}$. In fact $F^{Sp}$ must take vertices of $Q_{T,\epsilon}$ to vertices of $Q_{T,\epsilon'}$, edges to edges, faces to faces, etc. Indeed, $\mathfrak{s'} - \mathfrak{s} = h$ with $l(h) > 0$, we know that Since $F^{Sp}(\mathfrak{s'}) - F^{Sp}(\mathfrak{s}) = f_{\ast}(h)$, thus $F^{Sp}(\mathfrak{s'}) - F^{Sp}(\mathfrak{s}) = h'$ with $(f^{\ast})^{-1}(l) (h') > 0$. In particular, a supporting hyperplane for the set $Q_{T,\epsilon}$  must be mapped to a supporting hyperplane for $Q_{T,\epsilon'}$. Since vertices, edges and such are distigushed by the dimension of these supporting hyperplanes, we have the result we desired.\\
\ \\
\noindent To make use of the injectivity of $|f|$, note also that if $\mathfrak{s'} - \mathfrak{s} = k[h]$ where $k \in \Z$ and $[h]$ is primitive, then $F^{Sp}(\mathfrak{s'}) - F^{Sp}(\mathfrak{s}) = k f_{\ast}([h])$ as well. Thus, the divisibility of the difference of two $\mathrm{Spin}^{c}$ structures is preserved by $F^{Sp}$. \\
\ \\
\noindent  Suppose vertex $\mathfrak{s}$ for $Q_{T, \epsilon}$ is taken to vertex $\mathfrak{u}$ for $Q_{T,\epsilon'}$. Then each edge out of $\mathfrak{s}$ for $Q_{T,\epsilon}$ must be taken to an edge out of $\mathfrak{u}$. Furthermore, since the values of $|f(v)|$ are the divisibilities of the edges of the parallelepiped $Q_{T,\epsilon}$ out of $\mathfrak{s}$, when $|f|$ is injective each edge must be mapped to the unique edge emerging from $\mathfrak{u}$ with the same divisibility. Since the difference of the two vertices along this edge is a multiple of the primitive element $c_{v}$, we see that $c_{v} \rightarrow \pm c'_{v}$ under $f_{\ast}$ for every vertex $v$. As these span $H_{1}(S^{3}_{\Sigma_{T}})$, these described all the possibilities for $f_{\ast}$.\\
\ \\
\noindent However, by comparing $\overline{\theta}_{T,\epsilon}$ and $\overline{\theta}_{T,\epsilon'}$ we can see that none of these possibilities can occur for an orientation preserving diffeomorphism of pairs $f:(S^{3},\Sigma_{T,\epsilon}) \lra (S^{3},\Sigma_{T,\epsilon'})$. \\
\ \\
\noindent First, suppose there is an edge $e \in m$ with $\epsilon(e) = -\epsilon'(e)$. We may assume that $\epsilon(e) = +1$. Since $e \in m$, $e$ must be an edge directed from $b \in B$ to $w = m(b) \in W$. By our calculation in section 2, $\overline{\theta}_{T,\epsilon}(c_{b},c_{w}) = N_{-1, m}(\gamma) - N_{+1,\sim m}(\gamma)$. In our case $\gamma = e$, so $N_{-1,m}(\gamma) = 0$ while $N_{+1,\sim m}(\gamma) =0$. For $\epsilon'(e)  = -1$, we have  $N_{-1, m}(\gamma)=1$ while $N_{+1,\sim m}(\gamma) = 0$. Thus $\overline{\theta}_{T,\epsilon}(c_{b},c_{w}) = 0$ while $\overline{\theta}_{T,\epsilon'}(c'_{b},c'_{w})= 1$. From above the map $f$ can only map $c_{b} \rightarrow \pm c'_{b}$ and $c_{w} \rightarrow \pm c'_{w}$. As none of these choices can make the pairings equal, we see that there can be no map $f$ in this case.\\
\ \\
\noindent Second, assume that $\epsilon(e) = \epsilon'(e)$ for any edge $e \in m$. If $\epsilon \neq \epsilon'$ there is an edge $e \not\in m$ with $\epsilon(e) = - \epsilon'(e)$. Again, we assume that $\epsilon(e) = +1$. Since the edge $e$ is not in $m$, it must be directed from a vertex $w_{2} \in W$ to a vertex $b_{1} \in B$. Let $b_{2} = m^{-1}(w_{2})$ and $w_{1} = m(b_{1})$. We consider the directed path from $b_{2}$ to $w_{1}$:
$$
b_{2} \stackrel{g}{\Longrightarrow} w_{2} \stackrel{e}{\longrightarrow} b_{1} \stackrel{d}{\Longrightarrow} w_{1}
$$
By assumption $\epsilon(g) = \epsilon'(g)$ and $\epsilon(d) = \epsilon'(d)$. Once again we compute, $\overline{\theta}_{T,\epsilon}(c_{b},c_{w}) = N_{-1, m}(\gamma) - N_{+1,\sim m}(\gamma)$ $ = N - 1$ where $N$ is the number of $-1$'s in $\epsilon(g),\epsilon(d)$. On the other hand, $\overline{\theta}_{T,\epsilon'}(c'_{b},c'_{w}) = N'_{-1, m}(\gamma) - N'_{+1,\sim m}(\gamma)$ $= N - 0$ $=N$. Changing the signs of $c_{b}$ and/or $c_{w}$ can only change the value from $N-1$ to $1-N$. Thus, $f$ can only exist if $N$ satisfies either $N-1=N$, which has no solutions, or $1 - N = N$. The latter has no integer solutions, so we are done. $\Diamond$\\
\ \\
\noindent We note that using an $f$ assigning negative numbers to $b \in B$ and positive numbers to $w \in W$ makes $K_{T}$ have an alternating, prime, reduced diagram. By a well known theorem of Menasco, \cite{Men}, this means $K_{T}$ is also prime. Consequently, 

\begin{cor}
The prime, alternating knot $K_{T_{g}}$ defined by the matched tree $T_{g}$:
$$
w_{1} \Longleftarrow b_{1} \longleftarrow w_{2} \Leftarrow \cdots \leftarrow w_{2g} \Longleftarrow b_{2g}
$$
equipped with framing $f(w_{i}) = 2i+1$ and $f(b_{i}) = -2i$, bounds at least $2^{2g-1}$ inequivalent oriented surfaces of genus $g$. 
\end{cor}

\end{document}